\documentclass[11pt]{article}

\usepackage{amsmath}
\usepackage{amssymb}

 \newcommand{\beq}{\begin{equation}}
\newcommand{\eeq}{\end{equation}}

\newtheorem{theorem}{Theorem}[section]

\newtheorem{lemma}[theorem]{Lemma}
\newtheorem{corollary}[theorem]{Corollary}

\newcommand{\Cn}{{\mathbb  C}^n}

\newcommand{\Rnp}{{\overline{{\mathbb  R}_+^n}}}
\newcommand{\Zn}{{\mathbb  Z}_+^n}

\title{Special polyhedra for Reinhardt domains}
\author{Alexander Rashkovskii and Vyacheslav Zakharyuta}
\date{}

\begin{document}

\maketitle

\begin{abstract}
We show that every bounded hyperconvex Reinhardt domain can be approximated by special polynomial polyhedra defined by homogeneous polynomial mappings. This is achieved by means of approximation of the pluricomplex Green function of the domain with pole at the origin.
\end{abstract}

\section{Introduction and statement of results}

We consider here the following problem. {\sl Given a polynomially convex compact set $K\subset\Cn$ and its open neighborhood $U$, find a polynomial mapping $P=(P_1,\ldots,P_n)$ with $P^{-1}(0)\subseteq K$ and such that}
$$K\subset \{z\in \Cn:\: |P_k(z)|<1,\ 1\le k\le n\}\Subset U.$$
When $n=1$, its solution (D.~Hilbert, 1897) is known as the
Hilbert Lemniscate theorem; a proof can be found, e.g., in \cite{Ran}.
For $n=2$, the problem was solved in \cite{BLL} in the case of a circled set $K$ (that is, $\zeta K\subseteq K$ for any $\zeta$ from the closed unit disk), and the components of the mapping $P$ can be chosen as homogeneous polynomials of equal degree with unique common zero at the origin.

For arbitrary $n\ge 2$, it was shown in \cite{N}, \cite{N1} that weaker approximations are possible. Namely, either the mapping $P$ is allowed to have a small part of its zero set outside $K$ (for the general case of regular sets $K$), or with the number of components of the homogeneous mapping $P$ increased to $n+1$ (for the case of circled $K$). Even the approximation of the closed unit ball by the sublevel sets of $n$ homogeneous polynomials with $P^{-1}(0)=\{0\}$ was stated in \cite{N1} as an open problem.

\medskip
In this note, we solve the approximation problem by $n$ homogeneous polynomials for any multicircled (Reinhardt) polynomially convex compact set $K$. For pluripotential theory on multicircled sets and functions, see, for example, \cite{ARZ}, \cite{Za0}, \cite{Za}, \cite{Za03}.

\medskip
The above results in \cite{BLL}, \cite{N}, \cite{N1} make use of approximation of the pluricomplex Green function for a compact set with pole at infinity by logarithms of moduli of equidimensional polynomial mappings. Our approach is based on approximation of pluricomplex Green functions for bounded Reinhardt domains with pole at the origin. Note however that the Green function with logarithmic pole at $0$ for a circled domain $D$ can be extended by homogeneity $G(cz)=G(z)+\log|c|$ to the whole space, and the extension coincides on $\Cn\setminus D$ with the Green function for $\overline D$ with pole at infinity.

The main result is as follows.

\medskip
\begin{theorem}\label{theo:1} Let $G(z)$ be the pluricomplex Green function
of a bounded, polynomially convex Reinhardt domain
$D\subset\Cn$, with pole at $0$. Then there
exists a sequence of mappings $P^{(s)}:\Cn\to\Cn$
whose components $P_k^{(s)}$ are homogeneous polynomials of the same
degree $q_s$, with the only common zero at $0$, such that the sequence of
functions $v_s=q_s^{-1}\max_k\log|P_k^{(s)}|$ converges to $G$ uniformly
on $D\setminus\{0\}$, and the functions $v_s^+=\max\{v_s,0\}$ converge to the Green function for $\overline D$ with pole at infinity, uniformly on $\Cn$.
\end{theorem}

\medskip
As a consequence, we get an approximation of polynomially convex Reinhardt domains by special polyhedra defined by homogeneous polynomial mappings from $\Cn$ to $\Cn$, which solves, in particular, the problem posed in \cite[pp. 366-367]{N1}.

\medskip
\begin{corollary}\label{cor:1} For any bounded, polynomially convex
Reinhardt domain $D\subset\Cn$ and every $\epsilon>0$, there
exist $n$ homogeneous polynomials $p_k$ of the same degree, with the only common zero at $0$, such that
\begin{equation}\label{eq:Bn}\overline D\subset \{z:
|p_k(z)|<1,\ 1\le k\le n\}\subset (1+\epsilon)D.
\end{equation}
\end{corollary}

Finally, since any polynomially convex multicircled compact set is the intersection of domains satisfying the conditions of Corollary~\ref{cor:1}, we deduce

\medskip
\begin{corollary}\label{cor:2} For any polynomially convex
Reinhardt compact set $K\subset\Cn$ and every open neighborhood $U$ of $K$, there
exist $n$ homogeneous polynomials $p_k$ of the same degree, with the only common zero at $0$, such that
$K\subset \{z\in \Cn:\: |p_k(z)|<1,\ 1\le k\le n\}\Subset U$.
\end{corollary}

\section{Notation and preliminary results}

Let $D$ be a bounded hyperconvex domain. By $G_{a,D}$ we denote the {\it pluricomplex Green
function} of $D$ with pole at $a\in D$, that is, the upper envelope of negative plurisubharmonic functions $u$ in $D$ such that $ u(z)\le \log|z-a|+O(1)$ as $z\to a$.
The function $G$ is plurisubharmonic in $D$, continuous on
$\overline D$, maximal on $D\setminus\{a\}$, and
$G(z)=\log|z-a|+O(1)$ near $a$.

From now on, we specify $D$ to be a bounded, logarithmically convex and
complete (which amounts to being polynomially convex) Reinhardt ($n$-circled) domain in $\Cn$, and denote $G(z)=G_{0,D}(z)$.
Since the domain $D$ is $n$-circled, so is the function $G$.

Given $z\in\Cn$ and $\theta\in\Rnp$, denote
$S(\theta,z)=\sum_k\theta_k \log|z_k|$ and
let $h(\theta)=\sup\{S(\theta,z): z\in D\}$
be the characteristic function of the domain $D$ (the support function of the logarithmic image
$\log|D|$ of $D$).

\medskip
\begin{lemma}\label{prop:G} {\rm (cf., e.g., [9, Proposition 1.4.3], [1, Lemma 4])}
The Green function $G$ of a bounded
polynomially convex Reinhardt domain $D\subset\Cn$
with pole at $0$ has the representation
\begin{equation}\label{eq:repgreen}
G(z)=\sup\left\{S(\theta,z)-h(\theta):\:
\theta\in\Sigma\right\},\end{equation} where
$\Sigma=\{\theta\in\Rnp: \sum_k\theta_k=1\}$.
\end{lemma}

\medskip
{\it Proof.} Denote the right hand side of (\ref{eq:repgreen}) by $R(z)$. As is easy to see, it is plurisubharmonic in $D$, equal to zero on $\partial D$,
equivalent to $\log \left\vert
z\right\vert $ near $0$, and, since $R\left( cz\right) =
R\left( z\right) +\log \left\vert c\right\vert $, it is maximal on $D\setminus \left\{ 0\right\} $, so $R\left( z\right) \equiv G\left(
z\right) $ by the the Green function uniqueness property.
\medskip

\begin{lemma}\label{lem:Za03L6} For any $\epsilon>0$ and $t<0$, there
exist finitely many monomials $g_1,\ldots,g_m$ of the same degree $q$, such that
\begin{equation}\label{eq:locb}|G-v|<\epsilon\quad {\rm on\ }
D_{t,0}=\{z:\: t\le G(z)\le 0\},\end{equation} where
\begin{equation}\label{eq:maxg} v(z)=
q^{-1}\max\{\log|g_j(z)|:\: 1\le j\le m\},\end{equation} and the
maximum is attained for at most $n$ values of the indices $j$ at any
point $z$ on the level set $$\Gamma_t(v)=\{z\in D:\: v(z)=t\}.$$
\end{lemma}

\medskip

{\it Remark.} Lemma \ref{lem:Za03L6} can be deduced from a similar result \cite[Lemma~6]{Za03} on relative extremal functions. Note
that it was claimed there that, moreover, the maximum is attained
for at most $n$ functions $g_j$ on a set corresponding in our case to $\{ -1\le
v(z)\le 0\}$. However the proof of the claim has a gap,
and what is actually shown there is that the approximating function possesses this
property only on finitely many its level surfaces. That is why we
just follow the relevant arguments from the proof of \cite[Lemma~6]{Za03}.

\medskip

{\it Proof.} Representation
(\ref{eq:repgreen}), continuity of $G$, and compactness of $D_{t,0}$ imply the existence of
points $\theta^{(j)}\in\Sigma$, $1\le j\le m$, such that
$|G-u|<\epsilon/2$ on $D_{t,0}$, where
$$ u(z)=\max_{1\le j\le m} \{S(a^{(j)},z)-
b^{(j)}\},\quad
a^{(j)}=\theta^{(j)}, \quad b^{(j)}=h(\theta^{(j)}).$$
What we need to do is to
approximate the function $u$ by a similar one with rational
coefficients $\tilde a^{(j)}$ and to provide the required condition
about $n$ values. Consider the space $X\approx {\mathbb R}^{nm-p}\times {\mathbb
R}^{m}$ of matrices ${\mathcal M}=\left(a_{jk};\
b_j\right)\in{\mathbb R}^{nm}\times {\mathbb R}^{m}$ such that
$a_{jk}=0$ if $a^{(j)}_k=0$ ($p$ being the number of these
equations). For any $t<0$, there exists an algebraic set ${\mathcal A}_t\subset X$ such
that any system
$$ \sum_k a_{jk}x_k=b_j+t,\quad j\in J,\ |J|>n,$$ has no solution
for $(A,b)\in X\setminus{\mathcal A}_t$. Therefore, one can replace the
points $(a^{(j)},b^{(j)})$ by $(\tilde a^{(j)},\tilde b^{(j)})\in
X\setminus{\mathcal A}_t$ with rational $\tilde a^{(j)}\in r\Sigma$ for some rational $r$ close to $1$, such that the
function
$$ v(z)=\max_{1\le j\le m} \{S(\tilde a^{(j)},z)-
\tilde b^{(j)}\} $$ satisfies (\ref{eq:locb}) and the maximum is
attained for at most $n$ functions at any point of the set $\Gamma_t(v)$.

Finally, by choosing $N\in{\mathbb Z}_+$ so that $\tilde
a^{(j)}=N^{-1}k^{(j)}$, $rN\in\Zn$, and $k^{(j)}\in\Zn$ for all $j$, we get (\ref{eq:maxg}) with
monomials
$$g_j(z)=e^{-N\tilde b^{(j)}}z^{k^{(j)}}$$
of degree $q=rN$, which completes the proof. {\hfill$\square$\rm}

\medskip

The next point is a construction of precisely $n$ polynomials
approximating $G$ on $\Gamma_t(v)$. To this end, we use a
procedure from \cite[Lemma~2]{Za03}, see also \cite[Theorem~1]{AZa}.
Let $g_1,\ldots,g_m$ be the monomials from Lemma~\ref{lem:Za03L6}.
Given $s\in{\mathbb Z}_+$, let $g^{(s)}$ be a polynomial mapping
with the components
\begin{equation}\label{eq:qks}
g_k^{(s)}=\left(\sum_{j_1<\ldots<j_k}g_{j_1}^{s}\ldots
g_{j_k}^{s}\right)^{n!/k}, \quad k=1,\ldots, n.\end{equation}

\begin{lemma}\label{lem:Za03L2} {\rm \cite[Lemma~2]{Za03}} In the above notation, the sequence of functions
\begin{equation}\label{eq:vs}
v_s=(qsn!)^{-1}\max_{1\le k\le n}\log|g_k^{(s)}| \end{equation}
converges to $v$ uniformly on the level set $\Gamma_t(v)$.
\end{lemma}

\section{Proofs}

{\it Proof of Theorem~\ref{theo:1}}.
From Lemmas \ref{lem:Za03L6} and \ref{lem:Za03L2},  we derive
\begin{equation}\label{eq:collar}
|G(z)-v_s(z)|<\epsilon,\quad z\in \Gamma_{-1}(v)=\{v(z)=-1\},\ s\ge s_0(\epsilon).
\end{equation}

As follows from Lemma~\ref{prop:G}, the Green function $G$ satisfies
$ G(cz)=G(z)+\log|c|$ for all $c\in {\mathbb  C}$ such that
$cz\in D$. The functions $g_k^{(s)}$
defined by (\ref{eq:qks}) are homogeneous polynomials of degree
$q_s=qsn!$ and thus the function $v_s$ defined by (\ref{eq:vs}) has
the property
$ v_s(cz)=v_s(z)+\log|c|$ for all $c\in{\mathbb  C}$.  The homogeneity of both $G$ and $v_s$ extends (\ref{eq:collar}) to $\Cn\setminus\{0\}$ and implies the claimed convergence of the functions $v_s^+$.
{\hfill$\square$\rm}

\medskip
{\it Proof of Corollary \ref{cor:1}}. Take $\delta=\frac12\log(1+\epsilon)$. By Theorem~\ref{theo:1}, there exist $n$ homogeneous polynomials $P_k$ of degree $q$ such that
$|G(z)-q^{-1}\max_k\log|P_k(z)||<\delta$ for all $z\in \Cn$, $z\neq 0$. Then
\begin{equation}\label{eq:D1}\{z: |P_k|^{1/q}<e^{-\delta},\ 1\le k\le n\}\subset \overline D\subset \{z:
|P_k|^{1/q}<e^{\delta},\ 1\le k\le n\},
\end{equation}
which gives (\ref{eq:Bn}) with $p_k(z)=P_k(e^{-\delta}z)$.
{\hfill$\square$\rm}

Tek/Nat, University of Stavanger, 4036 Stavanger, Norway

{\sc E-mail}: alexander.rashkovskii@uis.no

\bigskip
Sabanci University, 34956 Tuzla, Istanbul, Turkey

{\sc E-mail}: zaha@sabanciuniv.edu


\begin{thebibliography}{11}

\bibitem{ARZ}
{\sc A. Aytuna, A. Rashkovskii and V. Zahariuta}, {\it Widths
asymptotics for a pair of Reinhardt domains}, Ann. Polon. Math. {\bf
78} (2002), 31-38.

\bibitem{AZa}
{\sc A. Aytuna and V. Zakharyuta}, {\it On Lelong--Bremermann lemma}, Proc. AMS {\bf 136} (2008), no. 5, 1733--1742.

\bibitem{BLL}
{\sc  T. Bloom, N. Levenberg, Yu. Lyubarskii}, {\it  A Hilbert
lemniscate theorem in ${\mathbb C}^2$}, Ann. Inst. Fourier (Grenoble) {\bf 58} (2008), no. 6, 2191--2220.

\bibitem{N}
{\sc S. Nivoche}, {\it Convexit\'e polynomiale, polyh\`edres
polynomiaux sp\'eciaux et applications}, C. R. Math. Acad. Sci.
Paris {\bf 342} (2006), no. 11, 825--830.

\bibitem{N1}
{\sc S. Nivoche}, {\it Polynomial convexity, special polynomial polyhedra and the pluricomplex Green function for a compact set in $\Cn$}, J. Math. Pures Appl. {\bf 91} (2009), 364--383.


\bibitem{Ran}
{\sc T. Ransford}, Potential Theory in the Complex Plane. Cmbbridge University Press, 1995.

\bibitem{Za0}
{\sc V.P. Zahariuta}, {\it Spaces of analytic functions and maximal
plurisubharmonic functions.} D.Sci. Dissertation, Rostov-on-Don,
1984.

\bibitem{Za} {\sc V. Zahariuta}, {\it Spaces of analytic functions and
Complex Potential Theory}, Linear Topological Spaces and Complex
Analysis {\bf 1} (1994), 74--146.

\bibitem{Za03} {\sc V. Zahariuta}, {\it On approximation by specialø
analytic polyhedral pairs}, Ann. Polon. Math. {\bf 80} (2003),
243--256.

\end{thebibliography}
\end{document}